\documentclass{elsarticle}

\usepackage{makeidx}
\usepackage{graphicx}

\usepackage{amsmath,amssymb}
\usepackage{appendix}
\usepackage[thref]{ntheorem}
\usepackage[numbers]{natbib}
\usepackage{xcolor}
\usepackage[font=small,labelfont={bf}]{caption}
\DeclareCaptionLabelSeparator{dot}{. } 
\usepackage{subfig}
\usepackage{csquotes}
\usepackage[normalem]{ulem}

\biboptions{sort&compress}

\newtheorem{rmk}{Remark}
\newtheorem{dfi}{Definition}
\newtheorem{thm}{Theorem}
\newtheorem{lem}{Lemma}

\newcommand{\bbR}{\mathbb R}
\newcommand{\RN}{\bbR^N}
\newcommand{\Rd}[1]{\bbR^{#1}}
\newcommand{\intlim}{\int}

\newcommand{\intt}{\intlim_0^t}
\newcommand{\intT}{\intlim_0^T}

\newcommand{\into}{\intlim_\Omega}
\newcommand{\intu}{\intlim_U}

\newcommand{\intdu}{\intlim_{\bndu}}

\newcommand{\intR}{ \intlim_{\mathbb{R}} }
\newcommand{\intRd}[1]{ \intlim_{\Rd{#1}} }
\newcommand{\intRN}{\intRd{N}}

\newcommand{\vek}[1]{\boldsymbol{#1}}
\newcommand{\nvek}{{\vek{\nu}}}
\newcommand{\bndo}{{\partial\Omega}}
\newcommand{\bndu}{{\partial U}}
\newcommand{\vx}{\vek{x}}
\newcommand{\llu}{ L^2(U) }
\newcommand{\llo}{ L^2(\Omega) }
\newcommand{\hu}{ H^1(U) }
\newcommand{\ho}{ H^1(\Omega) }
\newcommand{\lpThu}[1]{L^{ #1 }(0,T; \hu )}
\newcommand{\lpTho}[1]{L^{ #1 }(0,T; \ho )}

\newcommand{\lpu}  [1]{L^{ #1 }(U)}
\newcommand{\lpo}  [1]{L^{ #1 }(\Omega)}
\newcommand{\va}{\vek{a}}
\newcommand{\vb}{\vek{b}}
\newcommand{\D}{\mathrm{d}}
\newcommand{\dx}{\D\vx}
\newcommand{\dt}{\D t}

\newcommand{\ds}{\D s}

\newcommand{\nlpTho}[2]{\left\| #2 \right\|_{ \lpTho{#1} } }
\newcommand{\totdif}[2]{\frac{\D#1}{\D#2}}

\newcommand{\lpTlpo}[2]{L^{ #1 }(0,T; \lpo{#2} )}
\newcommand{\nlpTlpo}[3]{\left\| #3 \right\|_{ \lpTlpo{#1}{#2} } }
\newcommand{\nlpo} [2]{ \left\| #2 \right\|_{L^{ #1 }(\Omega)} }
\newcommand{\nlpThu}[2]{\left\| #2 \right\|_{ \lpThu{#1} } }

\begin{document}

\begin{frontmatter}
\title{Continuity in time of solutions of a phase-field model}
\author[sua]{Thomas G. Amler\corref{cor1}}
\ead{Thomas.Amler@kaust.edu.sa}
\author[mtu]{Nikolai~D.~Botkin}
\ead{botkin@ma.tum.de}
\author[mtu]{Karl-Heinz~Hoffmann}
\ead{hoffmann@ma.tum.de}
\author[sua]{Ibrahim Hoteit}
\ead{ibrahim.hoteit@kaust.edu.sa}

\cortext[cor1]{Corresponding author}

\address[sua]{Computer, Electrical and Mathematical Sciences and Engineering Division,
	4700 King Abdullah University of Science and Technology, Thuwal 23955-6900, Kingdom of Saudi Arabia.}

\address[mtu]{Technische Universit\"at M\"unchen, Zentrum Mathematik, Boltzmannstr.~3,
D-85748 Garching bei M\"unchen, Germany}

\begin{abstract}
	A phase field model proposed by G.\,Caginalp for the description of phase changes
	in materials is under consideration. It is assumed that the medium is located
	in a container with heat conductive walls that are not subjected
	to phase changes. Therefore, the temperature variable is defined both
	in the medium and wall regions, whereas the phase variable is only
	considered in the medium part. The case of Lipschitz domains in two and
	three dimensions is studied.
	We show that the temperature and phase variables are continuous
	in time functions with values in $L^2$ and $H^1$, respectively,
	provided that the initial values of them are  from $L^2$ and $H^1$, respectively.
	Moreover, continuous dependence of solutions on the initial data and boundary
	conditions is proved.
\end{abstract}
\end{frontmatter}

\section{Introduction}
\label{sec:intro}

Phase-field models have been introduced by G.\,Caginalp in \cite{Caginalp1986}
(see also e.g. \cite{Miranville2009} for generalizations)
as a relaxation of Stefan problems where the development of the sharp interface
between phases should be found exactly. In contrast to that, phase field models operate
with phase functions that assume values from -1 (solid state)
to 1 (liquid state) at each spatial point and change sharply
but smoothly over the solidification fronts so that the phase
interfaces become smoothed.
Mathematically speaking, such models consist of two coupled
parabolic equations describing  the temperature and phase fields
satisfying initial and boundary conditions.

Phase field models are frequently used to describe processes of melting,
solidification, evaporation, and condensation, which is directly related
to applications such as metal casting, design of cooling systems,
cryopreservation of living tissues etc.
(see e.g. \cite{EckHabil,BH2008Cryo}).
Phase field models are also appropriate for the description
of phase changes when modeling CO$_2$ sequestration.
It is assumed that the supercritical carbon dioxide,
CO$_2$  pressured to a phase between gas and liquid,
is injected into a saline
aquifer where it may either dissolve in the brine, react
with the dissolved minerals and the surrounding rocks, or
become trapped in the pore space of the aquifer.
In this relation, phase field models are useful
for the description of transitions between supercritical,
ordinary, and dissolved CO$_2$ phases.

In this paper, we study the case where a material subjected to a phase
change is located inside of a container with heat conductive walls that are not
subjected to phase changes. For example, this can be an ampoule filled with
a liquid containing living cells (see \cite[Sec. 2.2]{BH2008Cryo}) and subjected to
cooling applied to the outer surface but not immediately
to the liquid to be frozen.

The paper is structured as follows.
Section~\ref{sec:model} presents the mathematical model.
The definition of weak solutions and formulation of main results  are
given in Section~\ref{sec:main:res}.
The existence of weak solutions is shown in Section~\ref{sec:Galerkin}.
In Section \ref{sec:energy}, continuity in time of the phase function and
the validity of an energy equality are established.
The uniqueness of weak solutions and their continuous dependence on the
initial and boundary data are proved in Section~\ref{sec:uniq}.

\section{The model}
\label{sec:model}
First, the phase field model proposed by G.\,Caginalp in \cite{Caginalp1986}
will be recalled, then  generalizations considered in this paper will be introduced.

\subsection{Phase field model}
Let us outline the model proposed by G.\,Caginalp in \cite{Caginalp1986}.
Assume that a material
subjected to phase changes, e.g. liquid to
solid and back, occupies a region $\Omega\subset\RN$.
The evolution of the system is described in terms of the temperature
$u$ and a phase function $\phi$  satisfying  the following system of
non-linear parabolic equations:
\begin{equation}
	\label{model:o:eqs}
	\left.
	\begin{aligned}
		u_t
		+ \frac{l}{2}\,\phi_t
		- k \,\Delta u
		&
		= 0
		\\
		\tau \, \phi_t
		-2u
		+\frac{1}{2} \, \left(\phi^3-\phi\right)
		-\xi^2\Delta\phi
		&
		=0
	\end{aligned}
	\right\}
	\qquad \text{ in } \Omega\times(0,T).
\end{equation}
The first equation of \eqref{model:o:eqs} expresses the balance of heat energy.
Notice that this equation is scaled so that the term $u_t$ appears without any multiplier.
Here, $k$ is the scaled heat conductivity coefficient
which is assumed to be the same in the solid and liquid states, and
$l$ is the scaled specific latent heat of the phase change.
The second equation  of \eqref{model:o:eqs}
is derived from statistical physics, namely from the Landau-Ginzburg theory of
phase transitions \cite{Landau1958}.
In equilibrium, $\phi$ minimizes the following free energy functional:
\begin{equation*}
	F_u\{\phi\} =
	\intRN
	\left[
	\frac{\xi^2}2 |\nabla\phi|^2
	+ \frac{1}{8}
	\left(
	\phi^2 - 1
	\right)^2
	- 2u \, \phi
	\right],
\end{equation*}
where $\xi$ is a length scale, the thickness of the interfacial region
between liquid and solid.
In non-equilibrium, $\phi$ does not minimize $F_u\{\phi\}$ but
satisfies the gradient flow equation
$\tau\phi_t=-\delta F_u\{\phi\}$, where the symbol $``\delta"$ denotes
the variation of $F_u$ in $\phi$, which yields the second equation of \eqref{model:o:eqs}.
It is clear that the constant $\tau$ characterizes the relaxation time.

The following boundary conditions are usually imposed (see e.g. \cite{BH2008Cryo}):
\begin{equation}
	\label{model:o:BndC}
	\begin{aligned}
		-k \, \partial_\nvek u
		= \lambda(u-g),
		\qquad
		-\partial_\nvek\phi
		&
		= 0
		\quad
		\text{ on }\quad \bndo,
	\end{aligned}
\end{equation}
where $\nvek$ denotes the outer normal to $\Omega$,
and $\lambda$ is the scaled overall heat conductivity.

The specific of the case considered in this paper is that
the temperature $g$  is not directly applied to the boundary of the
liquid but to the outer surface of the container (ampoule).
To include the container into the model, still denote
the inner part occupied by the medium subjected to the phase change by
$\Omega$, and the solid walls of the container by $D$.
Moreover, let $U=D \cup \overline\Omega$ be the region occupied by the
whole system. In this case,
the phase function $\phi$ is defined in $\Omega$,
the temperature $u$ is defined in $U$,
and the boundary function $g$ is defined on $\bndu$.
Thus, the modification of the model \eqref{model:o:eqs} and \eqref{model:o:BndC}
reads as follows:
\begin{equation}
	\label{model:u:strong}
	\begin{aligned}
		\left.
		\begin{aligned}
			u_t
			+ \frac{l}{2}\,\phi_t
			- k_\Omega \,\Delta u
			&
			= 0
			\\
			\tau \, \phi_t
			-2u
			+\frac{1}{2} \, \left(\phi^3-\phi\right)
			-\xi^2\Delta\phi
			&
			=0
		\end{aligned}
		\right\}
		&  \qquad
		\text{ in } \Omega\times(0,T),
		\\
		\left.
		\begin{aligned}
			-k_D \, \partial_\nvek u_D
			&
			=
			-k_\Omega \, \partial_\nvek u_\Omega
			\\
			u_D
			&
			= u_\Omega
			\\
			-\partial_\nvek\phi
			&
			= 0
		\end{aligned}
		\right\}
		&  \qquad
		\text{ on }\bndo\times(0,T),
		\\[1ex]
		u_t
		- k_D \, \Delta u
		= 0
		& \qquad
		\text{ in } D\times(0,T),
		\\
		-k_D \, \partial_\nvek u
		= \lambda(u-g)
		&  \qquad
		\text{ on }\bndu\times(0,T),
		\\[1ex]
		u(\vx,0)=u^0(\vx) & \qquad \text{ in } U,
		\\[1ex]
		\phi(\vx,0)
		=\phi^0(\vx)
		&  \qquad \text{ in } \Omega.
	\end{aligned}
\end{equation}
Here the indices $D$ and $\Omega$ denote restrictions to the domains
$D$ and $\Omega$, respectively, and the symbol $\nvek$ is used to
indicate outward normals as well to $\Omega$ as to $U$, whenever it is not ambiguous.
Note that the matching conditions imposed on $u$ in \eqref{model:u:strong}
mean, in particular, the continuity of the temperature and the heat flux across $\bndo$.

The main result of this paper is that the problem \eqref{model:u:strong} admits a unique
weak solution $(u,\phi)$, and the mapping
\begin{equation*}
	\begin{split}
		&(u^0,\,\phi^0,\,g)
		\longrightarrow
		(u,\phi):\\
		&\llu \times \ho \times L^2(\bndu\times(0,T))
		\longrightarrow
		\mathcal C([0,T];\llu) \times \mathcal C([0,T];\ho)
	\end{split}
\end{equation*}
is continuous provided that $U$ and $\Omega$ are bounded Lipschitz domains in $\RN$, $N=2,3$.
A precise formulation of this result is given in the next section, see
Definition \ref{def:wsol} and Theorem \ref{thm:2reg}.

\section{Problem statement and main result}
\label{sec:main:res}

Problem \eqref{model:u:strong} will be studied
in a weak formulation. Denote $Q_T:=\Omega\times (0,T)$ and introduce the following spaces
according to  \cite{Magenes1} and \cite{wloka}:

\begin{equation}
	X_U:=\big\{
	u \in \lpThu 2
	\, : \;
	u_t \in L^2(0,T,(\hu)')
	\big\},
	\label{eq:def:X}
\end{equation}

\begin{equation}
	\begin{split}
		\mathcal C_s([0,T];\ho):=\big\{&
		\eta\in L^\infty(0,T,\ho):~~ t\rightarrow \left<\eta(t); \xi\right>\\
		&\mbox{is continuous on~~}  [0,T] \mbox{~~for each~~}
		\xi\in (\ho)^\prime
		\big\},
	\end{split}
	\label{eq:def:cs}
\end{equation}
where $\langle\cdot\,;\,\cdot\rangle$ denotes the dual pairing between  $\ho$
and $(\ho)^\prime$.

\begin{dfi}[Weak solutions]
	A pair $(u,\phi)$ of functions
	\begin{equation*}
		\begin{aligned}
			u \in X_U\qquad\text{and}\qquad
			\phi \in H^1(0,T;\llo) \cap \lpTho\infty
		\end{aligned}
	\end{equation*}
	satisfying the initial conditions
	\begin{equation*}
		u(0) = u^0
		\qquad\text{and}\qquad
		\phi(0) = \phi^0
	\end{equation*}
	is a weak solution of problem \eqref{model:u:strong},\; if the
	identities
	\begin{equation}\label{weak2}
		\begin{aligned}
			0
			&
			=
			\left< u_t ; \psi \right>_{X_U}
			+\intT\into \frac l 2 \,\phi_t \,  \psi
			+\intT\intu k\nabla u \cdot \nabla\psi
			+ \intT\int_{\partial U} \lambda\left(u-g\right)\psi,
			\\
			0
			&
			=
			\intT\into
			\left[
			\left(
			\tau \, \phi_t
			- 2u
			+ \frac{1}{2}\left(\phi^3-\phi\right)
			\right)
			\eta
			+ \xi^2 \, \nabla\phi \cdot \nabla\eta
			\right]
		\end{aligned}
	\end{equation}
	hold for all test functions $\psi\in L^2(0,T;H^1(U))$ and
	$\eta\in L^2(Q_T) \cap \lpTho 1$.
	\label{def:wsol}
\end{dfi}

\begin{rmk}
	\label{bemcs}
	Notice that the initial conditions in Definition \ref{def:wsol} have a sense
	because
	\begin{equation*}
		X_U \subset \mathcal C([0,T];\llu)
	\end{equation*}
	and
	\begin{equation*}
		H^1(0,T;\llo)\cap \lpTho\infty \subset \mathcal C_s([0,T];\ho),
	\end{equation*}
	see \cite[Th. 25.5]{wloka} for the first assertion
	and \cite[Chap. 3, Lemma 8.1]{Magenes1} for the second one.
\end{rmk}

The next theorem states the main result of this paper.

\begin{thm}[Main result]
	\label{thm:2reg}
	Let $U$ and $\Omega$ be bounded Lipschitz domains in $\RN, N=2,3,$
	with $\overline{\Omega} \subset U$, and $T>0$  be finite.
	If $u^0 \in \llu,$ $\phi^0\in \ho,$ and $g\in L^2(\bndu \times (0,T))$,
	then
	\begin{enumerate}
		\item
			There exists a unique weak solution $(u,\phi)$ of problem
			\eqref{model:u:strong} in the sense of Definition \ref{def:wsol}.
		\item
			The phase function $\phi$ has the additional regularity
			$\phi \in \mathcal C \left( [0,T];\ho \right)$.
		\item
			For all $s,\, t\in[0,T]$, the following energy equation holds:
			\begin{equation}
				\begin{aligned}
					&
					\into
					\left[
					\frac{1}{8}\phi(t)^4
					- \frac{1}{4}\phi(t)^2
					+ \frac{\xi^2}2|\nabla\phi(t)|^2
					\right]
					\\ & \quad
					=
					\into
					\left[
					\frac{1}{8}\phi(s)^4
					- \frac{1}{4}\phi(s)^2
					+ \frac{\xi^2}2|\nabla\phi(s)|^2
					\right]
					+\int_s^t \into
					\phi_t
					\left[
					2u- \tau\phi_t
					\right].
				\end{aligned}
				\label{eq:energy}
			\end{equation}
		\item
			Problem \eqref{model:u:strong} is well-posed in the sense that
			the mapping
			\begin{equation*}
				\begin{split}
					&(u^0,\,\phi^0,\,g) \longrightarrow (u,\phi):
					\\
					& \quad \llu \times \ho \times L^2(\bndu\times(0,T))
					\\
					& \quad \longrightarrow X_U \times \Big(H^1(0,T;\llo)\cap
					\mathcal C([0,T];\ho)\Big)
				\end{split}
			\end{equation*}
			is continuous.
	\end{enumerate}
\end{thm}

The proof of Theorem \ref{thm:2reg} is given in the following
sections.

\section{Existence of weak solutions}
\label{sec:Galerkin}

The existence of weak solutions of problem \eqref{model:u:strong}
can be proved in the similar way as Theorem 3.1.2 in \cite{EckHabil},
and therefore we only discuss changes in the construction of
approximate solutions and  a priori estimates.
For details related to the passage to the limit, we refer to \cite{EckHabil}.

The next lemma establishes some a priori estimates and the existence of weak solutions to problem
\eqref{model:u:strong}.

\begin{lem}
	\label{lem:exist}
	Let $U$ and $\Omega$ be bounded Lipschitz domains in $\RN, N=2,3$,
	with $\overline{\Omega} \subset U$, and $T>0$.
	If $u^0 \in \llu,$ $\phi^0\in \ho,$ and $g\in L^2(\bndu \times (0,T))$,
	then there exists a weak solution $(u,\phi)$ to problem \eqref{model:u:strong}
	in the sense of Definition \ref{def:wsol}.
\end{lem}

\noindent PROOF. Let $\{\zeta_i\}_{i=1}^\infty$ be a basis of $\hu$ which is orthonormal in
	$\lpu 2$, and $\{\omega_i\}_{i=1}^\infty$ a basis of $\ho$ which is orthonormal
	in $\llo$. Consider Galerkin approximations of the form
	\begin{equation}\label{approx2}
		u^{m}(\vx,t)=\sum\limits_{i=0}^m a_i^m(t) \, \zeta_i(\vx)
		\,,\qquad
		\phi^{m}(\vx,t)=\sum\limits_{i=0}^m b_i^m(t) \, \omega_i(\vx),
	\end{equation}
	where the functions $a_i^m(t)$ and $b_i^m(t)$ are to be determined.
	Let $\{g^m\} \subset$\newline $\mathcal C([0,T];L^2(\bndu))$ be a sequence
	such that $g^m \to g$ in $L^2(\bndu\times(0,T))$ as $m\to\infty$.

	To determine the functions $a_i^m$ and $b_i^m$, we require that $u^m$ and $\phi^m$
	satisfy the relations
	\begin{equation}
		\label{eq:weak:m}
		\begin{aligned}
			0
			&
			= \into \frac{l}{2} \, \phi^m_t(t) \, \psi^m
			+ \intu
			\left[
			u_t^m(t) \, \psi^m
			+ k \, \nabla u^m(t) \cdot \nabla\psi^m
			\right]
			\\
			&
			+ \int_{\partial U} \lambda
			\left(
			u^m(t)-g^m(t)
			\right)
			\psi^m,
			\\
			0
			&
			=
			\into
			\left[
			\tau \, \phi^m_t(t)
			- 2u^m(t)
			+ \frac{1}{2}\left( (\phi^m(t))^3-\phi^m(t)\right)
			\right]\eta^m
			\\
			&
			+
			\into
			\xi^2\nabla\phi^m(t) \cdot \nabla\eta^m
		\end{aligned}
	\end{equation}
	for	all test functions $\psi^m \in \mathrm{span}\{\zeta_j:j=1,\cdots,m\}$
	and $\eta^m \in \mathrm{span}\{\omega_j:j=1,\cdots,m\}$.

	Substituting ansatz \eqref{approx2} and
	test functions $\psi^m=\zeta_j$ and $\eta^m=\omega_j$,
	$j=1,\ldots,m$, into equations \eqref{eq:weak:m} yields the
	following system of ordinary differential equations for determining $a^m_j$ and $b^m_j$:
	\begin{equation}
		\label{coeff2}
		\begin{aligned}
			0
			&
			=
			\sum\limits_{i=1}^m
			\left\{
			\dot a_i^m(t)\intu \zeta_i \, \zeta_j
			+ k a_i^m(t)\intu\nabla\zeta_i \cdot \nabla\zeta_j
			+ \lambda\intdu a_i^m(t) \, \zeta_i \, \zeta_j
			\right\}
			\\
			&
			+
			\sum\limits_{i=1}^m
			\dot b_i^m(t)
			\, \frac{l}{2}
			\into \omega_i\,\zeta_j - \lambda\intdu g_m \zeta_j,
			\\
			0
			&
			=
			\sum\limits_{i=1}^m
			\left\{
			\tau \, \dot b_i^m(t)\into \omega_i \, \omega_j
			-\frac{1}{2}b_i^m(t)\into\omega_i \, \omega_j
			+\xi^2\,b_i^m(t)
			\into\nabla\omega_i \cdot \nabla\omega_j
			\right\}
			\\
			&
			-
			\sum\limits_{i=1}^m
			2a_i^m(t)\into\zeta_i \, \omega_j
			+\frac{1}{2}
			\into
			\left(
			\sum\limits_{i=1}^m
			b_i^m(t) \, \omega_i
			\right)^3
			\omega_j .
		\end{aligned}
	\end{equation}

	Since the set $\{\zeta_i\}$ is orthonormal in $L^2(U)$ and the set $\{\omega_i\}$ is orthonormal
	in $L^2(\Omega)$, equations \eqref{coeff2} can be rewritten as the following system of ODEs:
	\begin{equation}
		\label{ode2}
		\left[
		\begin{matrix}
			I_m & (l/2)M\\
			0 & \tau I_m
		\end{matrix}
		\right]
		\left[
		\begin{matrix}
			\dot{\vek{a}}^m(t)\\
			\dot{\vek{b}}^m(t)
		\end{matrix}
		\right]
		+\left[
		\begin{matrix}
			\vek{A}^m(\va^m(t),\,\vb^m(t))\\
			\vek{B}^m(\va^m(t),\,\vb^m(t))
		\end{matrix}
		\right]
		=\vek{0},
	\end{equation}
	where $M_{ij}:=\into \omega_i \, \zeta_j$, $i,j=1,\ldots,m$;
	$I_m$ denotes the $m\times m$ identity-matrix;
	and
	$\vek{A}^m\big(\va^m(t),\,\vb^m(t)\big)$ and $\vek{B}^m\big(\va^m(t),\,\vb^m(t)\big)$
	denote the terms of \eqref{coeff2} which do not comprise the time derivatives
	of the unknown functions.

	Initial conditions for $\va^m(t)$ and $\vb^m(t)$ are specified as follows:
	\begin{equation}\label{init2m}
		a_j^m(0) = \left( P^m_{\lpu2}u^0 \right)_j
		\qquad\text{ and }\qquad
		b_j^m(0) = \left( P^m_{H^1(\Omega)}\phi^0 \right)_j,
	\end{equation}
	where $P^m_{\lpu2}$ and $P^m_{\ho}$ denote the projectors onto the subspaces
	$\mathrm{span}\{\zeta_j:j=1,\cdots,m\}$
	and
	$\mathrm{span}\{\omega_j:j=1,\cdots,m\}$  in
	$\lpu2$ and $\ho$,  respectively.

	Notice that, for each fixed $m$, the functions $\vek A_j^m$ and
	$\vek B_j^m$ are analytic with respect to all their variables.
	Thus, the theory of ordinary differential equations provides that,
	for each $m$, there exist a nonempty time interval $[0,T_m]$ on
	which the initial value problem \eqref{ode2} and \eqref{init2m}
	admits a unique solution $[\va^m(t) , \vb^m(t)]$.

	To continue these solutions to any time interval, establish some
	independent on $m$ a priori estimates on $u^m$ and $\phi^m$.
	First, substitute $\psi^m=u^m(t)$ into the first equation of
	\eqref{eq:weak:m}, integrate it over $(0,t)$ for some $t\in(0,T_m]$,
	and use Young's inequality to obtain the estimate
	\begin{equation}
		\begin{aligned}
			&
			\frac 1 2 \intu
			|u^m(t)|^2
			+
			k
			\intt
			\intu
			|\nabla u^m|^2
			+
			\frac{\lambda} 2
			\intt \intdu |u^m|^2
			\\ & \quad
			\le
			\frac 1 2 \intu
			|u^0|^2
			+
			\frac{\lambda} 2
			\intt \intdu |g^m|^2
			+
			\epsilon
			\intt \into
			|\phi^m_t|^2
			+ \frac{l^2}{8\epsilon}
			\intt \intu
			|u^m|^2
		\end{aligned}
		\label{eq:estm:u}
	\end{equation}
	for $\epsilon > 0$.
	Next, substitute $\eta^m=\phi^m(t)$ into the second equation of
	\eqref{eq:weak:m} and proceed as before to obtain
	\begin{equation}
		\begin{aligned}
			&
			\frac \tau 2 \into |\phi^m(t)|^2
			+
			\intt \into
			\left[
			\frac{1}{2}
			\,
			|\phi^m|^4
			+
			\xi^2 |\nabla\phi^m|^2
			\right]
			\\ & \quad
			\le
			\frac \tau 2 \into |\phi^0|^2
			+
			\frac{1}{2} \intt \into |\phi^m|^2
			+
			2 \intt \intu |u^m|^2.
		\end{aligned}
		\label{eq:estm:phi}
	\end{equation}
	Now, substitute $\eta^m=\phi^m_t(t)$
	into the second equation of \eqref{eq:weak:m} and proceed as before
	to obtain
	\begin{equation}
		\begin{aligned}
			&
			(\tau - \epsilon)
			\intt \into |\phi^m_t(t)|^2
			+
			\into
			\left[
			\frac{1}{8}
			\,
			|\phi^m(t)|^4
			+
			\frac{\xi^2} 2 |\nabla\phi^m(t)|^2
			\right]
			\\ & \quad
			\le
			\into
			\left[
			\frac{1}{8}
			\,
			|\phi^0|^4
			+
			\frac{\xi^2} 2 |\nabla\phi^0|^2
			\right]
			+
			\frac 1 4 \into |\phi^m(t)|^2
			+
			\frac{1}{\epsilon} \intt \intu |u^m|^2
		\end{aligned}
		\label{eq:estm:phit}
	\end{equation}
	for $\epsilon > 0$.

	Choosing $\epsilon$ sufficiently small; multiplying inequalities
	\eqref{eq:estm:u}, \eqref{eq:estm:phi}, and \eqref{eq:estm:phit}
	by suitable constants; adding the resulting inequalities; and
	using the embedding $\ho \hookrightarrow \lpo 4$ yield the estimate
	\begin{equation}
		\begin{aligned}
			&
			\intu |u^m(t)|^2
			+
			\into
			\left[
			|\phi^m(t)|^2
			+
			|\phi^m(t)|^4
			+
			|\nabla\phi^m(t)|^2
			\right]
			\\ & \quad
			+
			\intt \intu |\nabla u^m|^2
			+
			\intt \intdu |u^m|^2
			+
			\intt \into |\phi^m_t(t)|^2
			\\ & \quad
			\le
			C
			\left[
			1
			+
			\intt \into
			|\phi^m|^2
			+
			\intt \intu
			|u^m|^2
			\right],
		\end{aligned}
		\label{eq:est:gron1}
	\end{equation}
	where the constant $C$ depends on $\tau$, $\xi$, $l$, $\lambda$,
	$\|u^0\|_{\lpu2}$, $\|\phi^0\|_{\ho}$, $\|g\|_{L^2(\bndu\times (0,T))}$,
	and $\epsilon$ but is independent on $m$.
	By Gronwall's inequality, the left-hand side of \eqref{eq:est:gron1}
	is bounded independently on $m$, which along with the embedding $\ho \hookrightarrow \lpo6$,
	$N\le 3$, implies the following assertions:
	\begin{equation}
		\label{bnd2}
		\begin{aligned}
			\{u^m\}
			&
			\text{ is bounded in }
			L^\infty\big(0,T_m;L^2(U)\big)\cap L^2\big(0,T_m;H^1(U)\big),
			\\
			\{u_t^m\}
			&
			\text{ is bounded in }
			L^2\big(0,T_m;(H^1(U))^\prime\big),
			\\
			\{\phi^m\}
			&
			\text{ is bounded in }
			L^\infty(0,T_m;H^1(\Omega)),
			\\
			\{\phi_t^m\}
			&
			\text{ is bounded in }
			L^2(0,T_m;L^2(\Omega)),
			\\
			\{(\phi^m)^3\}
			&
			\text{ is bounded in } L^\infty(0,T_m;L^2(\Omega)),
		\end{aligned}
	\end{equation}
	where the bounds are independent on $m$. This means that the approximate solutions
	can be continued to any interval $[0,T]$ keeping the above mentioned bounds. Now,
	the proof of the lemma can  be completed analogously to that of
	\cite[Th. 3.1.2]{EckHabil}. \quad $\square$

\section{Energy equalities}
\label{sec:energy}
In this section, two lemmas are proved.
Lemma~\ref{lem:equal:energy} states a slightly weaker than \eqref{eq:energy}
energy equality, which nevertheless implies the continuity in time of the
phase function: $\phi \in \mathcal C([0,T];\ho)$.
Lemma~\ref{lem:phi3t} completes the proof of the energy equality
\eqref{eq:energy}.
These results will be used in Section~\ref{sec:uniq} to show the claimed
uniqueness of weak solutions and their continuous dependence on the initial
and boundary data.

The proof of the next lemma is based on the techniques of
\cite[Chap 3, Sec 8.4, Lemma 8.3]{Magenes1}.

\begin{lem}[An energy equality]
	\label{lem:equal:energy}
	Let $(u,\phi)$ be a weak solution considered in Lemma~\ref{lem:exist}.
	Then, for all $s,\,t\in[0,T]$, the following energy equality holds:
	\begin{equation}
		\begin{aligned}
			&
			\into
			\left[
			|\nabla\phi(t)|^2
			-|\nabla\phi(s)|^2
			\right]
			=
			\frac2{\xi^2}
			\int_s^t\into
			\phi_t
			\left[
			2u
			- \tau\phi_t
			- \frac{1}{2}
			\left(\phi^3 - \phi \right)
			\right].
		\end{aligned}
		\label{eq:energy:weak}
	\end{equation}
	Moreover, equality \eqref{eq:energy:weak} implies that
	$\phi\in\mathcal{C}\left([0,T];H^1(\Omega)\right)$.
\end{lem}

\noindent PROOF. Assume that \eqref{eq:energy:weak} holds.
	Notice that the assertions of \eqref{bnd2} provide the integrability
	of the integrand in the right-hand side of \eqref{eq:energy:weak} and
	therefore the convergence of the integral to zero as $t\to s$,
	which proves the continuity of the function
	$t\rightarrow\into |\nabla\phi(t)|^2$. This, along with  the inclusion
	$\phi\in \mathcal C([0,T];\llo)$
	provided by the third and forth assertions of \eqref{bnd2}, implies
	the continuity of the function $t\rightarrow \|\phi(t)\|^2_{H^1(\Omega)}$.
	Let $t\in [0,T]$ be fixed, and $t_n\to t$ as $n\to \infty$. Denote
	$\epsilon_n=\|\phi(t_n)-\phi(t)\|^2_{H^1(\Omega)}$ and compute
	$\epsilon_n = \|\phi(t_n)\|^2_{H^1(\Omega)} + \|\phi(t)\|^2_{H^1(\Omega)} -
	2 \big(\phi(t_n),\phi(t)\big)_{H^1(\Omega)}$.
	Utilizing that $\phi\in \mathcal C_s([0,T];H^1(\Omega))$
	(see \eqref{eq:def:cs} and remark~\ref{bemcs}) implies the convergence $\epsilon_n\to 0$, which means
	that  $\phi(t_n)\to \phi(t)$ in $H^1(\Omega)$.

	Now go on to the proof of equality \eqref{eq:energy:weak} and introduce some notation.
	Similar to \cite{Magenes1}, denote  the inner product
	in $\llo$ or the dual pairing in $\ho'\times \ho$ by $\langle\cdot;\cdot\rangle$,
	and denote the inner product in $L^2(\mathbb{R};\llo)$ or the dual pairing
	in $L^2(\mathbb{R};\ho')\times L^2(\mathbb{R};\ho)$ by $(\cdot;\cdot)$.
	Assume that $\phi(t)$ is defined on $\mathbb{R}$ and possesses the same properties
	as $\phi(t)$ on $[0,T]$. This can be achieved trough continuation of $\phi(t)$ by means of
	reflections. In the following, set $s=0$ and $t=t_0$ in equation \eqref{eq:energy:weak}.

	Define a bilinear form $a:\ho\times\ho\to\bbR$
	and an operator $A:\ho\to\ho'$ by
	\begin{equation}
		a(v,\,w)=
		\left<
		Av \,;\, w
		\right>
		= \into \, \nabla v \cdot \nabla w \, \dx,
		\qquad
		v,\,w \in \ho.
		\label{eq:A:def}
	\end{equation}
	Then the second equation of \eqref{weak2} is equivalent to the relation
	\begin{equation}
		\intT
		\left<
		A \phi(t)\, ; \, \eta(t)
		\right>
		\dt
		=
		\frac1{\xi^2}
		\intT\into
		\left[
		2u
		- \tau \phi_t
		- \frac 1 2 \left( \phi^3 -\phi \right)
		\right]
		\eta.
		\label{eq:weak:A}
	\end{equation}

	For $t_0\in(0,T]$ and $\delta \in (0,t_0/2)$, define
	the function $Q_\delta:\bbR\to[0,1]$ by
	\begin{equation}\label{def:Q}
		Q_\delta(t):=
		\begin{cases}
			1 \quad
			&
			\mathrm{for }\,\, t\in[\delta, t_0-\delta]
			\\
			0 \quad
			&
			\mathrm{for }\,\, t\notin [0,t_0]
			\\
			\text{linear on}
			&
			[0,\delta]\cup[t_0-\delta,t_0]\,\,
			\text{ and continuous on } \mathbb{R}.
		\end{cases}
	\end{equation}

	Let $\{\rho_n\}$ be a sequence of non-negative even regularizing
	functions with $\intR\rho_n(t) \dt = 1$ and
	$\mathrm{supp} \rho_n \subset\subset (-1/n, 1/n)$. Denote $\sigma_n=\rho_n*\rho_n$.

	Moreover, as it was already mentioned (see remark~\ref{bemcs}),
	\begin{equation}
		\phi \in \mathcal C_s([0,T];\ho).
		\label{eq:phi:scalcont}
	\end{equation}

	Since $A$ is a symmetric operator, it holds
	\begin{equation}\label{intgradphit}
		\begin{aligned}
			0
			&
			= \intR \totdif {} t
			\left<
			A (\rho_n \ast Q_\delta\,\phi)
			\, ; \,
			\rho_n \ast Q_\delta\,\phi
			\right> \dt
			\\ &
			= 2 \left(
			A (\rho_n \ast Q_\delta\,\phi)
			\, ; \,
			\rho_n \ast Q_\delta\,\phi_t
			\right)
			+
			2
			\left(
			A(\rho_n \ast Q_\delta\,\phi)
			\, ; \,
			\rho_n \ast Q_\delta'\,\phi
			\right)
			\\ &
			=
			2
			\left(
			A (\rho_n \ast Q_\delta\,\phi)
			\, ; \,
			\rho_n \ast Q_\delta\,\phi_t
			\right)
			+
			2
			\left(
			\rho_n \ast A (Q_\delta - Q_0) \phi
			\, ; \,
			\rho_n \ast Q_\delta'\,\phi
			\right)
			\\ &
			+
			2
			\left(
			\rho_n \ast (A \, Q_0\,\phi)
			\, ; \,
			\rho_n \ast Q_\delta'\,\phi
			\right).
		\end{aligned}
	\end{equation}
    Notice that $A\phi\in L^2(Q_T)$ because of the second equation of \eqref{model:u:strong}
    and the assertions of \eqref{bnd2}. This implies that
	\begin{equation*}
		A\left( \rho_n\ast Q_\delta\,\phi \right) \in L^2(\bbR;\llo),
	\end{equation*}
	since $\rho_n$ and $Q_\delta$ depend only on time, and $A$ is time
	independent. Therefore, all duality brackets on the right-hand
    side of equation \eqref{intgradphit} present the inner product
    of $L^2(\bbR;\llo)$.

	Consider the passage to the limit
    in each term of the right-hand side of \eqref{intgradphit}
    as $\delta \to 0$ and then as $n\to\infty$. The objective is to show
    the following three relations:
	\begin{align}
		&2 \left(
		A (\rho_n \ast Q_\delta \phi)
		\, ; \,
		\rho_n \ast Q_\delta \phi_t
		\right)
		\longrightarrow
		- \frac 2 {\xi^2}
		\intT\into
		\phi_t
		\left[
		\tau \phi_t
		-2u
		+ \frac 1 2 \left( \phi^3 -\phi \right)
		\right],
		\label{term1lim}
		\\[2ex]
		&\left|
		\left(
		\rho_n \ast A (Q_\delta - Q_0) \phi
		\, ; \,
		\rho_n \ast Q_\delta ' \phi
		\right)
		\right|
		\longrightarrow 0,
		\label{term2lim}
		\\[2ex]
		&2
		\left(
		\rho_n \ast (A \, Q_0 \, \phi)
		\, ; \,
		\rho_n \ast Q_\delta' \phi
		\right)
		\longrightarrow
		\left<
		A\, \phi(0)
		\, ; \,
		\phi(0)
		\right>
		-
		\left<
		A\, \phi(t_0)
		\, ; \,
		\phi(t_0)
		\right>
		,
		\label{term3lim}
	\end{align}
	as first $\delta \to 0$, and then $n\to\infty$, which
	along with equation \eqref{intgradphit} completes the proof of the lemma.

	\emph{Proof of \eqref{term1lim}}:
	Consider the limit as $\delta \rightarrow 0$.
	By properties of convolutions, it holds
	\begin{equation}
		A\left( \rho_n\ast (Q_\delta \, \phi) \right)
		= \rho_n\ast\left(Q_\delta \, A \, \phi\right)
		\longrightarrow
		\rho_n\ast\left(Q_0 \, A \, \phi\right)
		\qquad\text{ in }L^2(\bbR;\llo)
		\label{eq:Aphi:l2}
	\end{equation}
	because these functions have compact supports in $t$.
	Moreover, $\rho_n \ast (Q_\delta \, \phi_t) \rightarrow \rho_n \ast (Q_0 \, \phi_t)$
	in $L^2(\bbR;\llo)$, and therefore, the first term of the right-hand side of the last
	equation of \eqref{intgradphit} satisfies
	\begin{equation}
		\label{eq:lim:1}
		2 \left(
		A (\rho_n \ast Q_\delta\,\phi)
		\, ; \,
		\rho_n \ast Q_\delta\,\phi_t
		\right)
		\longrightarrow
		2 \left(
		A ( \rho_n \ast Q_0\,\phi)
		\, ; \,
		(\rho_n \ast Q_0\,\phi_t)
		\right)
	\end{equation}
	as $\delta \to 0$.
	Relation \eqref{eq:lim:1}, the second equation of \eqref{weak2},
	and time independence of $A$ imply that
	\begin{equation}
		\label{intterm1}
		\begin{aligned}
			&
			2 \left(
			A (\rho_n \ast Q_\delta\,\phi)
			\, ; \,
			\rho_n \ast Q_\delta\,\phi_t
			\right)
			\\
			& \quad
			\longrightarrow
			- \frac 2 {\xi^2}
			\intT\into
			(\sigma_n \ast (Q_0\,\phi_t))
			\left[
			\tau \phi_t
			-2u
			+ \frac 1 2 \left( \phi^3 -\phi \right)
			\right]
		\end{aligned}
	\end{equation}
	as $\delta \to 0$.
	Moreover, note that
	\begin{equation}
		2 \sigma_n \ast (Q_0 \, \phi_t)(\cdot) =
		\int_{\mathbb R^+} 2\sigma_n(\cdot -s) \, Q_0(s) \, \phi_t(s)
		\longrightarrow Q_0 \, \phi_t
		\label{eq:lim:1:sigma}
	\end{equation}
	in $L^2(\bbR;\llo)$ as $n\to\infty$,
	because $\int_{\mathbb R^+} \sigma_n(s)\, \ds = 1/2$.
	Therefore, relations \eqref{intterm1} and \eqref{eq:lim:1:sigma} yield
	\eqref{term1lim}.

	\emph{Proof of \eqref{term2lim}}:
	Using H\"older's inequality yields
	\begin{equation}
		\label{eq:term2:1:hoelder}
		\begin{aligned}
			&
			\left|
			\left(
			\rho_n \ast A (Q_\delta - Q_0) \phi
			\, ; \,
			\rho_n \ast Q_\delta'\,\phi
			\right)
			\right|
			\\
			& \quad
			\le
			\left\|
			\rho_n \ast A (Q_\delta - Q_0) \phi
			\right\|_{L^\infty(\bbR;\llo)}
			\cdot
			\left\|
			\rho_n \ast Q_\delta'\,\phi
			\right\|_{L^1(\bbR;\llo)}.
		\end{aligned}
	\end{equation}
	Using again the relation $A\phi\in L^2(Q_T)$ and the properties of the convolution
    operator, one can show that
	the first factor on the right-hand side of estimate \eqref{eq:term2:1:hoelder}
	tends to zero.
	The second factor on the
	right-hand side of \eqref{eq:term2:1:hoelder} can be estimated using Young's inequality as follows:

	\begin{equation*}
		\begin{aligned}
			\left\| \rho_n \ast Q_\delta' \, \phi \right\|_{L^1(\bbR;\llo)}
			& \le
			\left\| \rho_n \right\|_{L^1(\bbR)}
			\cdot
			\left\| Q_\delta' \, \phi \right\|_{L^1(\bbR;\llo)}
			\\ & \le
			\left\| Q_\delta' \right\|_{L^1(\bbR)}
			\cdot
			\left\| \phi \right\|_{L^\infty(0,T;\llo)}.
		\end{aligned}
	\end{equation*}
	The right-hand side is bounded, because $\|\rho_n\|_{L^1(\bbR)} =1$, $\|Q_\delta'\|_{L^1(\bbR)} =2$, and
	$\mathrm{supp}\,Q_\delta \subset[0,T]$.
	Thus, the right-hand side of  \eqref{eq:term2:1:hoelder} tends to zero as
	$\delta\to0$, and relation \eqref{term2lim} is proved.

	\emph{Proof of \eqref{term3lim}}:
	Notice that $\rho_n\ast(A\phi) \in \mathcal C([0,T];\llo)$
	because $A\phi \in L^2(Q_T)$.
	Thus, the function
	\begin{equation*}
		t \mapsto \left< \left( \sigma_n \ast A\,Q_0\,\phi
		\right) (t) \, ; \, \phi(t) \right>
	\end{equation*}
	is continuous since $\phi\in\mathcal C([0,T];\llo)$.
	Accounting for the evenness of $\rho_n$, the left-hand side of
	\eqref{term2lim} admits the following passage to the limit
	as $\delta\to0$:
	\begin{equation}\label{lim0T}
		\begin{aligned}
			&
			2 \left( \rho_n \ast (A\,Q_0\,\phi) \, ; \,
			\rho_n \ast Q_\delta'\,\phi \right)
			\\ & \quad
			\longrightarrow
			2 \left< (\sigma_n\ast(A\,Q_0\,\phi))(0)
			\,;\, \phi(0) \right>
			-2 \left< (\sigma_n\ast(A\,Q_0\,\phi))(t_0)
			\,;\, \phi(t_0) \right>.
		\end{aligned}
	\end{equation}
	It remains to show that
	\begin{equation}
		\begin{aligned}
			2
			\left<
			(\sigma_n\ast(A\,Q_0\,\phi))(t_0)
			\,;\,
			\phi(t_0)
			\right>
			&
			\longrightarrow
			\left<
			A\,\phi(t_0) \,;\, \phi(t_0)
			\right>,
			\\
			2
			\left<
			(\sigma_n\ast(A\,Q_0\,\phi))(0)
			\,;\,
			\phi(0)
			\right>
			&
			\longrightarrow
			\left<
			A \, \phi(0)
			\,;\,
			\phi(0)
			\right>
		\end{aligned}
		\label{eq:term:aim}
	\end{equation}
	as $n\to\infty$.
	Consider the first relation of \eqref{eq:term:aim}.
	By the definition of $Q_\delta$ and $\sigma_n$,
	we have
	\begin{equation*}
		\begin{aligned}
			2 \left< (\sigma_n\ast(A\,Q_0\,\phi))(t_0)
			\,;\, \phi(t_0) \right>
			= 2 \int^{t_0}_0 \sigma_n(t)
			\left< (A\,Q_0\,\phi)(t_0-t) \,;\,
			\phi(t_0) \right> \dt
		\end{aligned}
	\end{equation*}
	and $\int^{t_0}_0 \sigma_n(t) \, \dt =\frac{1}{2}$
	because $\sigma_n$ is even.
	From the definition of the operator  $A$ (see \eqref{eq:A:def}),
	it follows that
	\begin{equation}
		\label{eq:term:laplace}
		\begin{aligned}
			&
			2
			\left< (\sigma_n\ast(A\,Q_0\,\phi))(t_0) \,;\,
			\phi(t_0) \right>
			-
			\left< A \, \phi(t_0) \,;\, \phi(t_0) \right>
			\\
			&
			\quad
			=
			2 \int^{t_0}_0 \sigma_n(t)
			\left< A \, \phi(t_0-t) - A \, \phi(t_0)
			\,;\, \phi(t_0) \right>
			\, \dt
			\\
			& \quad
			= 2 \int^{t_0}_0 \sigma_n(t) \,
			a(\phi(t_0-t) - \phi(t_0),\, \phi(t_0)) \, \dt
			\\
			& \quad
			\longrightarrow 0\quad \mbox{ as }\quad n\to\infty,
		\end{aligned}
	\end{equation}
	because the function $t\mapsto a(\phi(t_0-t),\, \phi(t_0))=
	\left<\phi(t_0-t),\, A\phi(t_0)\right>$ is continuous,
	see \eqref{eq:phi:scalcont}.
	Thus, Lemma \ref{lem:equal:energy} is proved. \quad $\square$

The next lemma proves an integral form of the relation
$(\phi^3-\phi)\phi_t = (\phi^4/4 -\phi^2/2)_t$, the chain rule.

\begin{lem}
	\label{lem:phi3t}
	If $\phi \in H^1(0,T;\llo)\cap  \mathcal C([0,T];\ho)$,
	then
	\begin{equation}
		\int_s^t \into
		\phi_t
		\left[
		\phi^3 - \phi
		\right]
		=
		\frac{1}{4}
		\into
		\left[
		|\phi(t)|^4
		- |\phi(s)|^4
		\right]
		-
		\frac{1}{2}
		\into
		\left[
		|\phi(t)|^2
		- |\phi(s)|^2
		\right]
		\label{eq:phi:poly}
	\end{equation}
	for all $s, \, t\in[0,T]$.
	\label{lem:poly:chain}	
\end{lem}

\noindent PROOF. Consider the case $s=0$.
	Since $\phi\in H^1(0,T;\llo)$, it holds
	\begin{equation}
		\intt \into
		\phi_t
		\,
		\phi
		=
		\frac 1 2
		\into
		\left[
		\phi(t)^2
		-
		\phi(0)^2
		\right],
		\label{eq:phi2:aim}
	\end{equation}
	see e.g. \cite[Th. 1.67]{novotny}.

	To show the relation
	\begin{equation}
		\intt \into
		\phi_t
		\,
		\phi^3
		=
		\frac 1 4
		\into
		\left[
		\phi(t)^4 - \phi(0)^4
		\right],
		\label{eq:phi4:aim}
	\end{equation}
	notice that $H^1(0,T;\llo)\cap \lpTho\infty \hookrightarrow \mathcal C([0,T];\lpo4)$,
	which follows from \cite[Sec. 8, Cor 4]{jsimon} and the compact embedding
	$\ho\hookrightarrow\hookrightarrow\lpo4$ for $N\le 3$.
	Let $\{\phi_\epsilon\}$ be a sequence of smooth functions such that
	\begin{equation}
		\phi_\epsilon \longrightarrow \phi
		\qquad
		\text{in } H^1(0,T;\llo)\cap \mathcal C([0,T];\ho)\quad \text{as } \epsilon \to 0.
		\label{eq:phi4:approx}
	\end{equation}
	The chain rule yields
	\begin{equation}
		\intt \into
		(\phi_\epsilon)_t
		\,
		\phi_\epsilon^3
		=
		\frac 1 4
		\into
		\left[
		\phi_\epsilon(t)^4 - \phi_\epsilon(0)^4
		\right].
		\label{eq:phi4:eps}
	\end{equation}
	The passage to the limit on the both sides of equation \eqref{eq:phi4:eps}, as $\epsilon\to0$,
	has to be done.
	The left-hand side is being processed as follows:
	\begin{equation}
		\begin{aligned}
			&
			\intt \into
			\left[
			\phi_\epsilon^3 \, (\phi_\epsilon)_t
			-
			\phi^3 \, \phi_t
			\right]
			=
			\intt \into
			\left[
			\left(
			\phi_\epsilon^3
			-
			\phi^3
			\right)\,
			(\phi_\epsilon)_t
			+
			\phi^3
			\left(
			(\phi_\epsilon)_t
			-
			\phi_t
			\right)
			\right].
		\end{aligned}
		\label{eq:phi4:diff}
	\end{equation}
	The first summand in the integral on the right-hand side can be estimated using
	H\"older's inequality
	and the formula $ \phi_\epsilon^3 - \phi^3 = (\phi_\epsilon - \phi) \zeta$
	with $\zeta=\phi_\epsilon^2 + \phi^2 + \phi_\epsilon \, \phi$.
	It holds
	\begin{equation}
		\begin{aligned}
			&
			\left|
			\intt\into (\phi_\epsilon - \phi) \, \zeta\cdot \,
			(\phi_\epsilon)_t
			\right|^2
			\le
			\intt\into
			(\phi_\epsilon - \phi)^2 \, \zeta^2
			\cdot
			\intt\into
			\left|
			(\phi_\epsilon)_t
			\right|^2
			\\
			& \quad
			\le
			\intt
			\left\{
			\left[
			\into
			\left|
			\phi_\epsilon - \phi
			\right|^6
			\right]^{1/3}
			\cdot
			\left[
			\into
			\left|
			\zeta^3
			\right|
			\right]^{3/2}
			\right\}
			\dt
			\cdot
			\intt\into
			\left|
			(\phi_\epsilon)_t
			\right|^2
			\\ & \quad
			\le
			C
			\nlpTho \infty { \phi_\epsilon - \phi }^2
			\cdot
			\nlpTlpo {9/2} 3 \zeta ^{9/2}
			\cdot
			\|(\phi_\epsilon)_t\|^2_{L^2(Q_T)}.
		\end{aligned}
		\label{eq:phi4:term1}
	\end{equation}
	It is not hard to prove tat $\zeta \in \lpTlpo \infty {3}$. Really,
	Young's inequality and the embedding $\ho\hookrightarrow\lpo6$
	yield the estimate
	\begin{equation}
		\begin{aligned}
			\into \zeta^3
			&
			= \into
			\left[
			\phi_\epsilon^2
			+ \phi_\epsilon \, \phi
			+ \phi^2
			\right]^3
			\\ &
			\le
			\into
			\phi_\epsilon^6
			+
			\into
			\phi^6
			+
			C
			\sum_{k=1}^5
			\into
			\left|
			\phi_\epsilon^k
			\,
			\phi^{6-k}
			\right|
			\\ &
			\le
			\into
			\phi_\epsilon^6
			+
			\into
			\phi^6
			+
			C
			\sum_{k=1}^5
			\left[
			\frac{k}{6}
			\into
			\phi_\epsilon^6
			+
			\frac{6-k}{6}
			\into
			\phi^{6}
			\right]
			\\ &
			\le C
			\left[
			\nlpTho\infty \phi^{6}
			+
			\nlpTho\infty {\phi_\epsilon}^{6}
			\right].
		\end{aligned}
		\label{eq:zeta:l3}
	\end{equation}
	Thus, estimates \eqref{eq:phi4:term1} and \eqref{eq:zeta:l3}
	yield
	\begin{equation}
		\intt\into
		(\phi_\epsilon - \phi)
		\,
		\zeta
		\,
		(\phi_\epsilon)_t
		\longrightarrow 0
		\qquad
		\text{ for } \epsilon\to0.
		\label{eq:phi4:term1:fin}
	\end{equation}

	Consider the second summand in the integral on the right-hand side of equation \eqref{eq:phi4:diff}.
	Processing it similarly to the first term yields
	\begin{equation}
		\begin{aligned}
			&
			\left|
			\intt \into
			\phi^3
			\big(
			(\phi_\epsilon)_t
			-
			\phi_t
			\big)
			\right|^2
			\le
			\intt\into \phi^6
			\cdot
			\intt\into
			\left|
			(\phi_\epsilon)_t
			-
			\phi_t
			\right|^2
			\\ & \quad
			\le
			C \, T \,
			\nlpTho \infty \phi^6
			\cdot
			\nlpTlpo 2 2 { (\phi_\epsilon)_t - \phi_t }^2
			\\
			& \quad
			\longrightarrow 0
			\qquad \text{ as } \epsilon\to 0.
		\end{aligned}
		\label{eq:phi4:term2}
	\end{equation}

	Now, estimates \eqref{eq:phi4:term1:fin} and \eqref{eq:phi4:term2}
	along with equation \eqref{eq:phi4:diff} yield
	\begin{equation}
		\intt \into
		\left[
		\phi_\epsilon^3 \, (\phi_\epsilon)_t
		-
		\phi^3 \, \phi_t
		\right]
		\longrightarrow 0
		\qquad
		\text{as } \epsilon\to0.
		\label{eq:phi4:diff:lim}
	\end{equation}

	Using the approximation property \eqref{eq:phi4:approx},
	and the embedding $\ho\hookrightarrow\lpo6$ yields the estimate
	\begin{equation}
		\begin{aligned}
			&
			\into
			\left[
			\phi^4(t) -\phi_\epsilon(t)^4
			\right]
			=
			\into
			\left[
			\left(
			\phi(t) - \phi_\epsilon(t)
			\right)
			\zeta(t)
			\right]
			\\ & \quad
			\le
			\nlpo2{\phi(t) - \phi_\epsilon(t)}
			\cdot
			\nlpo 2{\zeta(t)}
			\\ & \quad
			\longrightarrow 0
		\end{aligned}
		\label{eq:phi4:diff4:lim}
	\end{equation}
	for all $t\in[0,T]$ as $\epsilon\to0$.
	Here we have used the abbreviation $\zeta = \phi^3 + \phi^2\phi_\epsilon + \phi^2\phi_\epsilon
	+ \phi_\epsilon^3 \in \mathcal C([0,T];\llo)$.

	Thus, equation \eqref{eq:phi4:aim} follows from equation \eqref{eq:phi4:eps} and
	relations \eqref{eq:phi4:diff:lim} and \eqref{eq:phi4:diff4:lim}, which proves Lemma~\ref{lem:phi3t}.
\quad $\square$

Finally, the energy equality \eqref{eq:energy} follows from Lemmas
\ref{lem:equal:energy} and \ref{lem:poly:chain}.

\section{Uniqueness and continuous dependency on the data}
\label{sec:uniq}

The next lemma completes the proof of Theorem \ref{thm:2reg}.

\begin{lem}
	Let $U$ and $\Omega$ be bounded Lipschitz domains in $\RN, N=2,3,$
	with $\overline{\Omega} \subset U$, and $T>0$.
	Let $u^0_i \in L^2(U)$, $\phi^0_i \in \ho$, and
	$g_i \in L^2(\bndu\times(0,T)), \, i=1,2,$
	be given initial and boundary data of problem \eqref{model:u:strong};
	$(u_i,\phi_i)$ weak solutions of problem \eqref{model:u:strong} in
	the sense of Definition \ref{def:wsol}. If
	$\bar u=u_1 -u_2$ and $\bar\phi = \phi_1 -\phi_2$,
	then
	\begin{equation}
		\left.
		\begin{aligned}
			&
			\|\bar u\|_{\mathcal C([0,T];\llu)}
			\\
			&
			\nlpThu 2 {\bar u}
			\\
			&
			\|\bar \phi\|_{\mathcal C([0,T];\ho)}
			\\
			&
			\|\bar \phi\|_{H^1(0,T;\llo)}
		\end{aligned}
		\right\}
		\le
		F\big( \|\bar u^0\|_{L^2(U)}, \,
		\|\bar\phi^0\|_{\ho}, \,
		\|\bar g\|_{L^2(\bndu\times(0,T))}\big),
		\label{eq:cont}
	\end{equation}
	where $F:[0,\infty)^3 \to [0,\infty)$ is a continuous function
	with $F(0,0,0)=0$.
	\label{lem:uniq}
\end{lem}

\noindent PROOF. Denote $\zeta = \phi^2_1 + \phi_1\,\phi_2 + \phi^2_2 $ and bear in mind that
	$\zeta$ is non-negative.
	Obviously, $\bar u$ and $\bar\phi$ satisfy the equations
	\begin{equation}\label{weakbar2}
		\begin{aligned}
			0
			&
			=
			\left< \bar u_t ; \psi \right>_{X_U}
			+ \frac{l}{2}\intT\into \bar\phi_t \, \psi
			+\intT\intu k \, \nabla \bar u \cdot \nabla\psi
			+ \intT\int_{\partial U} \lambda\left(\bar u- \bar g\right)\psi
			\\
			0
			&
			=
			\intT\into\left[
			\left(
			\tau \bar\phi_t
			+
			- 2\bar u
			+ \frac{\zeta}{2} \bar \phi
			- \frac{1}{2}\, \bar\phi
			\right)
			\eta
			+ \xi^2\,\nabla\bar\phi \cdot \nabla\eta\right]
		\end{aligned}
	\end{equation}
	for all test functions $\psi\in L^2(0,T;H^1(U))$ and $\eta\in \lpTho 2$.
	In particular, test functions $\psi=\chi_{(0,t)} \, \bar u$,
	$\eta=\chi_{(0,t)} \, \bar\phi$, and $\eta=\chi_{(0,t)} \, \bar\phi_t$
	will be considered to obtain three estimates.
	Here, $\chi_{(0,t)}$ denotes the characteristic function of the
	interval $(0,t)$.
	Notice that the proofs of Lemmas~\ref{lem:equal:energy} and
	\ref{lem:poly:chain} show that $\eta=\chi_{(0,t)} \, \bar\phi_t$ is an
	admissible test function.

	Substituting $\psi = \chi_{(0,t)} \, \bar u$ into the first equation
	of \eqref{weakbar2} and using Young's inequality yield the estimate
	\begin{equation}\label{baru2}
		\begin{aligned}
			&
			\frac 1 2 \intu |\bar u(t)|^2
			+ \intt\intu k |\nabla\bar u|^2
			+ \frac \lambda 2 \intt\intdu|\bar u|^2
			\\ & \quad
			\le
			\frac 1 2 \intu |\bar u^0|^2
			+ \frac \lambda 2 \intt\intdu |\bar g|^2
			+ \intt\into
			\left[
			\epsilon |\bar\phi_t|^2
			+ \frac{l^2}{16\,\epsilon} |\bar u|^2
			\right]
		\end{aligned}
	\end{equation}
	for any $\epsilon > 0$.
	Substituting $\eta=\chi_{(0,t)} \, \bar\phi$ into the second equation
	of \eqref{weakbar2} and applying Young's inequality yield
	\begin{equation}\label{barphi2}
		\begin{aligned}
			&
			\frac \tau 2 \into |\bar\phi(t)|^2
			+ \intt\into
			\left[
			\frac \zeta 2 |\bar\phi|^2
			+ \xi^2 |\nabla\bar\phi|^2
			\right]
			\\ & \quad
			\le \frac \tau 2 \into |\bar\phi^0|^2
			+ \intt\into\left[ |\bar u|^2 + \frac 3 2 |\bar\phi|^2 \right].
		\end{aligned}
	\end{equation}
	Substituting $\eta=\chi_{(0,t)} \, \bar\phi_t$ into the second equation
	of \eqref{weakbar2} and applying Young's inequality yield
	\begin{equation}\label{barenergy}
		\begin{aligned}
			&
			\into\left[
			\frac{1}{4}\left|\bar\phi^0\right|^2
			+\frac{\xi^2}2|\nabla\bar\phi(t)|^2
			\right]
			+
			\int_0^{t}\into
			\frac{\tau}{4}
			\,
			|\bar\phi_t|^2
			\\
			& \quad
			\le
			\into
			\left[
			\frac{1}{4}|\bar\phi(t)|^2
			+\frac{\xi^2}2|\nabla\bar\phi^0|^2
			\right]
			+
			\intt\into
			\left[
			\frac 2 \tau \, \bar u^2
			+ \frac 1 {4\tau} \, \zeta^2\, \bar\phi^2
			\right]
			.
		\end{aligned}
	\end{equation}
	Using  H\"older's inequality  and the embedding $\ho\hookrightarrow\lpo6$ yield
	\begin{equation}
		\begin{aligned}
			&
			\intt\into
			\left[
			\zeta^2\, \bar\phi^2
			\right]
			\le
			\intt
			\left\{
			\left[
			\into |\zeta|^3
			\right]^{2/3}
			\cdot
			\left[
			\into |\bar\phi|^6
			\right]^{3}
			\right\}
			\\ & \quad
			\le
			\nlpTlpo \infty 3 \zeta^2
			\intt
			\nlpo 6 {\bar\phi(t)}^2
			\\ & \quad
			\le
			C
			\nlpTlpo \infty 3 \zeta^2
			\cdot
			\intt \into
			\left[
			|\bar\phi(t)|^2
			+
			|\nabla\bar\phi(t)|^2
			\right].
		\end{aligned}
		\label{eq:zeta:bar}
	\end{equation}
	Now, choose $\epsilon \in (0,\tau/4)$ in \eqref{baru2} and combine estimates
	\eqref{baru2}, \eqref{barphi2}, \eqref{barenergy}, and \eqref{eq:zeta:bar}
	to obtain
	\begin{equation}
		\begin{aligned}
			&
			\intu |\bar u(t)|^2
			+ \into |\nabla\bar\phi(t)|^2
			+ \intt\intu |\nabla\bar u|^2
			+ \intt\into |\bar\phi_t|^2
			+ \intt\intdu|\bar u|^2
			\\ & \quad
			\le
			C \into |\bar\phi(t)|^2
			+
			C
			\left[
			\intu |\bar u^0|^2
			+ \into |\nabla\bar\phi^0|^2
			+ \intt\intdu |\bar g|^2
			\right]
			\\ & \quad
			+
			C
			\left[
			\intt\intu |\bar u|^2
			+ C_\zeta^2
			\intt \into
			\left[
			|\bar\phi|^2
			+
			|\nabla\bar\phi|^2
			\right]
			\right],
		\end{aligned}
		\label{eq:bar:est1}
	\end{equation}
	where $C_\zeta:= \nlpTlpo \infty 3 \zeta$, and the constant $C$ depends on
	$k,\,l,\,\lambda,\xi,$ and $\tau$ only.
	To eliminate the first term on the right-hand side of \eqref{eq:bar:est1},
	multiply inequality \eqref{barphi2} by a constant greater than $2C/\tau$
	and add the resulting inequality to \eqref{eq:bar:est1} (remember that $\zeta$ is non-negative).
	This yields
	\begin{equation}
		\begin{aligned}
			&
			\intu |\bar u(t)|^2
			+ \into
			\left[
			|\bar\phi(t)|^2 + |\nabla\bar\phi(t)|^2
			\right]
			\\ & \quad
			+ \intt\intu |\nabla\bar u|^2
			+ \intt\into |\bar\phi_t|^2
			+ \intt\intdu|\bar u|^2
			\\ & \quad
			\le
			C
			\left\{
			C_0
			+
			(1+C_\zeta^2)
			\intt
			\left(
			\intu |\bar u|^2
			+
			\into
			\left[
			|\bar\phi|^2 + |\nabla\bar\phi|^2
			\right]
			\right)
			\right\},
		\end{aligned}
		\label{eq:bar:est2}
	\end{equation}
	where $C$ is sufficiently large, and
	\begin{equation}
		C_0 :=
		\intu |\bar u^0|^2
		+ \into
		\left[
		|\bar\phi^0|^2+|\nabla\bar\phi^0|^2
		\right]
		+ \intt\intdu |\bar g|^2.
		\label{eq:def:C0}
	\end{equation}
	Applying Gronwall's inequality to \eqref{eq:bar:est2} yields
	\begin{equation}
		\begin{aligned}
			&
			\intu |\bar u(t)|^2
			+ \into
			\left[
			|\bar\phi(t)|^2
			+
			|\nabla\bar\phi(t)|^2
			\right]
			\\ & \quad
			+ \intt\intu |\nabla\bar u|^2
			+ \intt\into |\bar\phi_t|^2
			+ \intt\intdu|\bar u|^2
			\\ & \quad
			\le
			C_0 \, C
			\left[
			1+ T \, C (1+C_\zeta^2)
			\exp \left( T \, C (1+C_\zeta^2) \right)
			\right].
		\end{aligned}
		\label{eq:bar:est3}
	\end{equation}
	Due to the a priori estimates derived in Section \ref{sec:Galerkin}
	(see the assertions of \eqref{bnd2}), the constant $C_\zeta$ is bounded in
	terms of the problem data only.
	Thus, estimate \eqref{eq:bar:est3} proves the continuous dependence
	of weak solutions on the problem data.
	In particular,
	inequality \eqref{eq:bar:est3} and the definition of $C_0$
	(see \eqref{eq:def:C0}) imply uniqueness of weak solutions.\quad $\square$

\section{Simulation}

The following simulation deals with a plastic ampoule used for freezing
small tissue samples (see \cite{BH2008Cryo}).
The diameter of the ampoule is equal to 1\,cm, the height
equals 5\,cm, and the wall thickness equals 0.1\,cm. The ampoule is filled
with water, and is being cooled with the rate of 1$^\circ$C/s applied to the outer surface
of the ampoule. The mass and thermal characteristics of plastic and water
are easily available in reference books. The length scale $\xi$ and the
relaxation time $\tau$ are equal to $0.03$ and $0.005$, respectively.

\begin{center}
	\begin{figure}[h!]
		\setlength{\unitlength}{1cm}
		\begin{picture}(9,5)
			\put(0,0){
			\put(0.2,0){\includegraphics[height=50mm]{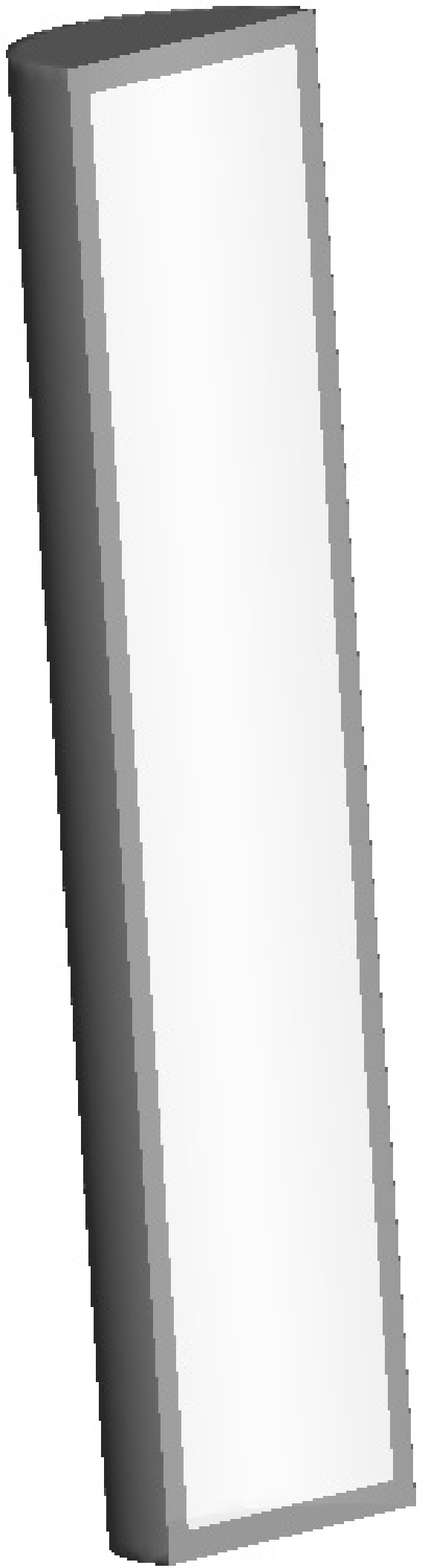}}
			\put(3.0,0){\includegraphics[height=50mm]{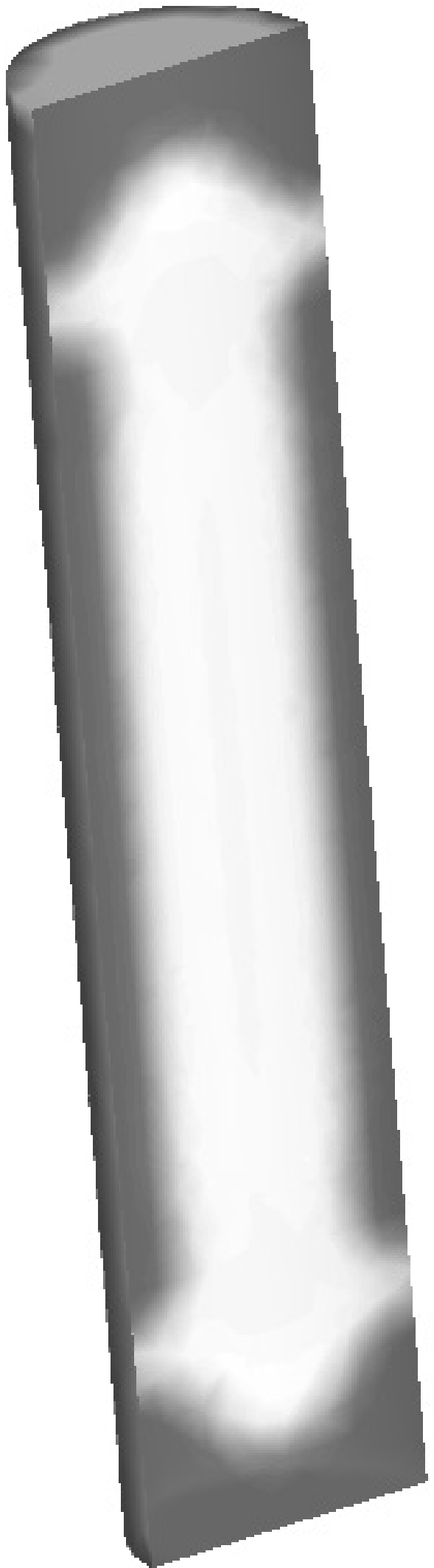}}
			\put(5.3,0){\includegraphics[height=50mm]{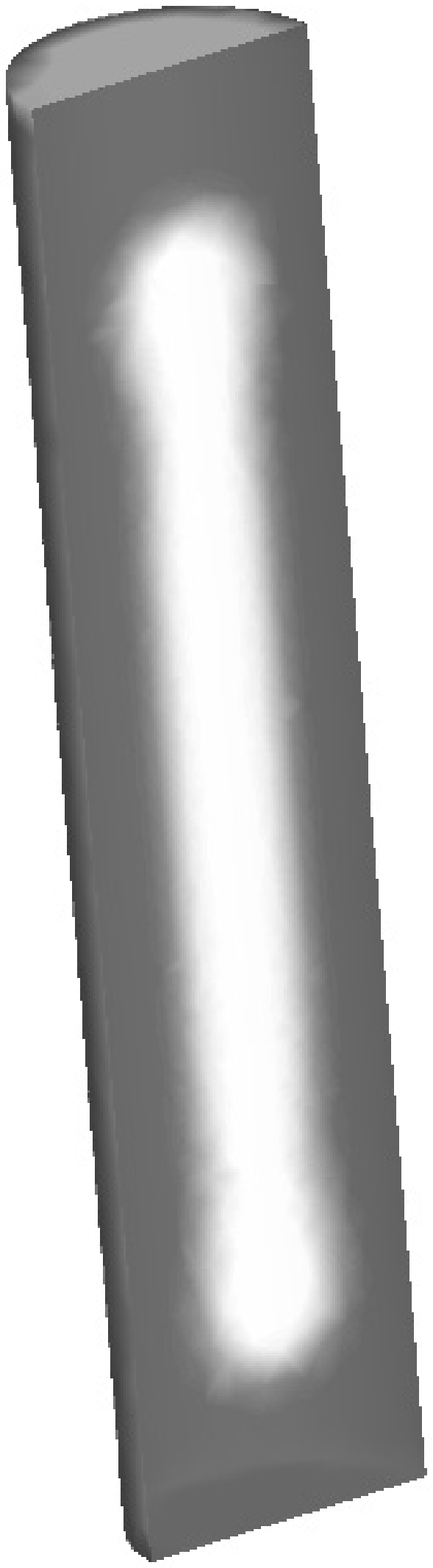}}
			\put(7.9,0){\includegraphics[height=50mm]{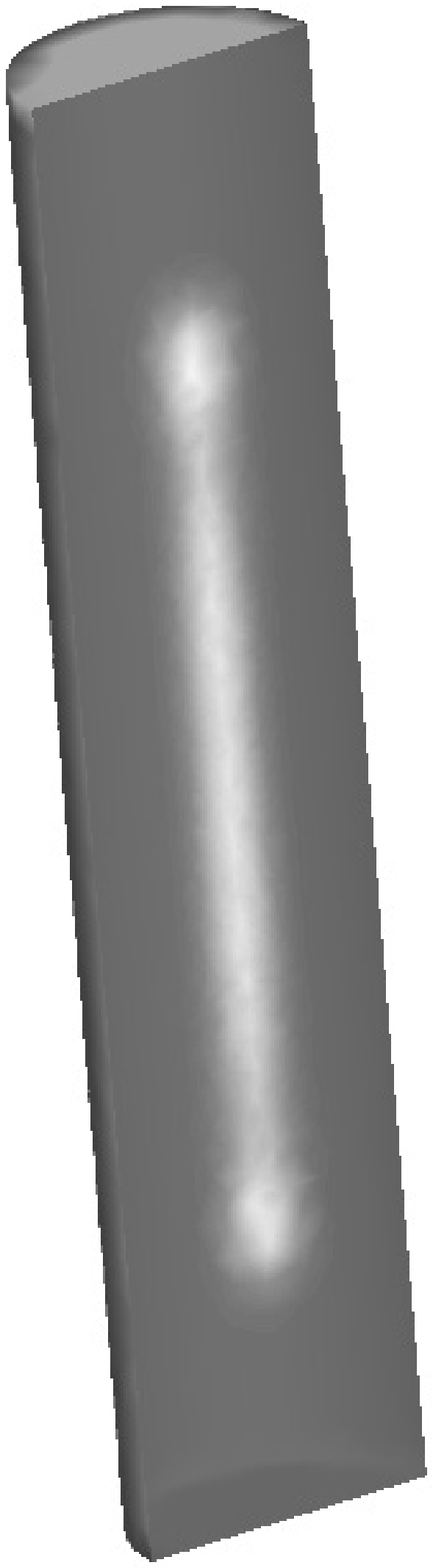}}
			\put(10.5,0){\includegraphics[height=50mm]{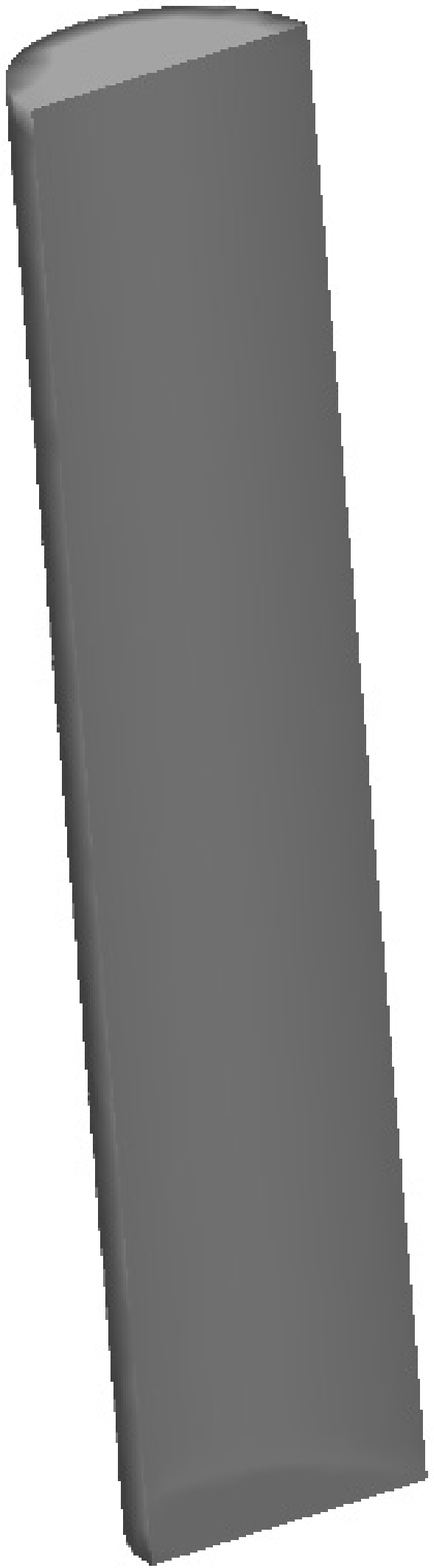}}

			\put(0,3){$t_0$}
			\put(2.8,3){$t_1$}
			\put(5.2,3){$t_2$}
			\put(7.6,3){$t_3$}
			\put(10.3,3){$t_4$}
			}
		\end{picture}
		\caption{Time development of the unfrozen (light) and frozen (dark) parts.
		Notice that $\phi\approx 1$ in the light region, and  $\phi\approx -1$ in the dark one.
		There is a small transition zone between the unfrozen and frozen parts.}
		\label{fig1}
	\end{figure}
\end{center}

The first picture of Figure~\ref{fig1} shows an axial cut
of the whole ampoule (region $U$).  The others show the water part
(region $\Omega$) only. The unfrozen part is shown in light, whereas the
frozen one becomes dark.

\section{Conclusion}

We have considered a phase field model describing phase changes of
a medium located in a container with heat conductive walls that are
free of phase changes. It is shown that the temperature and the phase variables are
continuous functions with values in $L^2(U)$ and $H^1(\Omega)$, respectively,
provided that the initial data are  from $L^2(U)$ and $H^1(\Omega)$,
respectively.
Moreover, solutions depend continuously on the initial  and
boundary data of the problem.\\

\noindent {\bf Acknowledgements.}
The work was motivated by the Project SPP 1253 of the German Research Society (DFG),
and was supported by Award No. KSA-C0069/UK-C0020,
made by King Abdullah University of Science and Technology (KAUST).


\newpage

\begin{thebibliography}{}

	\bibitem{Caginalp1986}
		G. Caginalp. An analysis of a phase field model of a free boundary.
		{\em Arch. Rat. Mech. Anal.} 92(3):205-245, 1986.

	\bibitem{Miranville2009}
		A. Miranville and R. Quintanilla. Some generalizations of the Caginalp phase-field system.
		{\em Applicable Analysis}, 88(6):877-894, 2009.


	\bibitem{EckHabil}
		C. Eck. {\em A Two-Scale Phase Field Model for Liquid-Solid Phase Transitions of Binary Mixtures
		with Dendritic Microstructrue.} Habilitation, Universit\"at Erlangen-N\"urnberg, July 2004.

	\bibitem{BH2008Cryo}
		K.-H. Hoffmann and N. D. Botkin. Optimal control in cryopreservation of cells and tissues.
		{\em Advances in Mathematical Sciences and Applications}, (29):177-200, 2008.

	\bibitem{Landau1958}
		L. D. Landau and E. M. Lifshitz. {\em Statistical Physics}. Reading, MA: Addison-Wesley, 1958.

	\bibitem{Magenes1}
		J. L. Lions and E. Magenes. {\em Probl\'emes aux limites non homog\`enes et applications},
		volume~1. Dunod, Paris, 1968.

	\bibitem{wloka}
		J. Wloka. {\em Partielle Differentialgleichungen -- Sobolevr\"aume und Randwertaufgaben}. B.G.
		Teubner, Stuttgart, 1982.

	\bibitem{novotny}
		A. Novotny and I. Straskraba. {\em Introduction to the Mathematical Theory of Compressible
		Flow}, volume 27 of {\em Oxford Lecture Series in Mathematics and Its Applications}. Oxford
		University Press, August 2004.

	\bibitem{jsimon}
		J. Simon. Compact sets in the space ${L}^p(0,{T}; {B})$. {\em Annali di Matematica Pura ed Applicata},
		146(1):65 - 96, 2005.


\end{thebibliography}
\end{document}